\documentclass[10pt]{article}
\usepackage{cite}
\usepackage{mathrsfs}
\usepackage{amsfonts}
\usepackage{amsmath}
\usepackage{amsfonts,amssymb}
\usepackage{dsfont}
\usepackage{curves}
\usepackage{mathrsfs}
\usepackage{pifont}
\usepackage{amssymb}
\allowdisplaybreaks

\numberwithin{equation}{section}

\date{}

\def\BigRoman{\uppercase\expandafter{\romannumeral\number\count 255 }}
\def\Romannumeral{\afterassignment\BigRoman\count255=}

\begin{document}
\title{Spectral radius and $k$-factor-critical graphs
\thanks{This work was supported by the Natural Science Foundation of Shandong Province, China (ZR2023MA078)}
}
\author{\small  Sizhong Zhou$^{1}$\footnote{E-mail address: zhousizhong@just.edu.cn (S. Zhou)},
Zhiren Sun$^{2}$\footnote{E-mail address: 05119@njnu.edu.cn (Z. Sun)}, Yuli Zhang$^{3}$\footnote{Corresponding
author. E-mail address: zhangyuli\_djtu@126.com (Y. Zhang)}\\
\small  $1$. School of Science, Jiangsu University of Science and Technology,\\
\small  Zhenjiang, Jiangsu 212100, China\\
\small  $2$. School of Mathematical Sciences, Nanjing Normal University,\\
\small  Nanjing, Jiangsu 210023, China\\
\small  $3$. School of Science, Dalian Jiaotong University,\\
\small  Dalian, Liaoning 116028, China\\
}

\maketitle
\begin{abstract}
\noindent For a nonnegative integer $k$, a graph $G$ is said to be $k$-factor-critical if $G-Q$ admits a perfect matching for any $Q\subseteq V(G)$
with $|Q|=k$. In this article, we prove spectral radius conditions for the existence of $k$-factor-critical graphs. Our result generalises one previous
result on perfect matchings of graphs. Furthermore, we claim that the bounds on spectral radius in Theorem 3.1 are sharp.
\\
\begin{flushleft}

(2020) Mathematics Subject Classification: 05C50, 05C70
\end{flushleft}
\end{abstract}

\section{Introduction}

Let $G$ be a graph with vertex set $V(G)$ and edge set $E(G)$ which has neither multiple edges nor loops. We denote by $n=|V(G)|$ the order of
$G$. The number of odd components in $G$ is denoted by $o(G)$. The number of connected components in $G$ is denoted by $\omega(G)$. For any
$v\in V(G)$, we denote by $d_G(v)$ the degree of $v$ in $G$. For any $D\subseteq V(G)$, we denote by $G[D]$ the subgraph of $G$ induced by $D$,
and by $G-D$ the graph formed from $G$ by removing the vertices in $D$ and their incident edges. A vertex subset $X$ of $G$ is called an
independent set if any two members of $X$ are not adjacent in $G$. For a graph $G$ of order $n$, the adjacency matrix $A(G)$ of $G$ is the
$n$-by-$n$ matrix in which entry $a_{ij}$ is 1 or 0 according to whether $v_i$ and $v_j$ are adjacent or not, where $V(G)=\{v_1,v_2,\cdots,v_n\}$.
The eigenvalues of the adjacency matrix $A(G)$ are also called the eigenvalues of $G$. The largest eigenvalue of $G$, denoted by $\rho(G)$, is
called the spectral radius of $G$. For an integer $k\geq2$, the sequential join $G_1\vee G_2\vee\cdots\vee G_k$ of graphs $G_1,G_2,\cdots,G_k$
is the graph with vertex set $V(G_1\vee G_2\vee\cdots\vee G_k)=V(G_1)\cup V(G_2)\cup\cdots\cup V(G_k)$ and edge set $E(G_1\vee G_2\vee\cdots\vee
G_k)=\{e: e\in E(G_i) \ \mbox{for some} \ 1\leq i\leq k \ \mbox{or an unordered pair between} \ V(G_i) \ \mbox{and} \ V(G_{i+1}) \ \mbox{for some} \
1\leq i\leq k-1\}$. Let $K_n$ and $P_n$ denote the complete graph and the path of order $n$, respectively.

Let $a$ and $b$ be two integers with $0\leq a\leq b$. Then a spanning subgraph $F$ of a graph $G$ is called an $[a,b]$-factor if
$a\leq d_F(v)\leq b$ for all $v\in V(G)$. When $a=b=1$, an $[a,b]$-factor is simply called a 1-factor (or a perfect matching). For a nonnegative
integer $k$, a graph $G$ is said to be $k$-factor-critical if $G-Q$ admits a 1-factor for any $Q\subseteq V(G)$ with $|Q|=k$.

Egawa and Furuya \cite{EF} studied the existence of perfect matchings in star-free graphs. Brouwer and Haemers \cite{BH}, O \cite{Os} presented
eigenvalue conditions for graphs to possess perfect matchings. Johansson \cite{J} verified an El-Zah\'ar type condition for the existence of
$[1,2]$-factors in graphs. Zhou, Bian and Pan \cite{ZBP}, Zhou, Sun and Liu \cite{ZSLo}, Zhou, Sun and Yang \cite{ZSY}, Zhou \cite{Zd,Zs,Zp},
Gao and Wang \cite{GW}, Gao, Chen and Wang \cite{GCW}, Liu \cite{L}, Wang and Zhang \cite{WZd,WZi}, Wu \cite{Wp}, Li and Miao \cite{LM}
demonstrated some theorems on $[1,2]$-factors in graphs. Wang and Zhang \cite{WZr} derived a sufficient condition for the existence of
factor-critical graphs. Lv \cite{Lv}, Zhou \cite{Za1}, Zhou, Pan and Xu \cite{ZPX}, Zhou, Liu and Xu \cite{ZLX}, Wu \cite{Wa} posed some sufficient
conditions for graphs to be factor-critical graphs. Gu and Liu \cite{GL}
established a relationship between eigenvalues and factor-critical graphs. Ananchuen and Plummer \cite{AP} showed a result for the existence of
3-factor-critical graphs. Enomoto, Plummer and Saito \cite{EPS} investigated a relationship between neighborhoods of independent sets and
$k$-factor-critical graphs. Plummer and Saito \cite{PS} showed some characterizations for graphs to be $k$-factor-critical graphs. Wang and Yu
\cite{WY} proved that a $k$-connected 3-$\gamma$-edge-critical claw-free graph $G$ with minimum degree at least $k+1$ and $k\equiv|V(G)|$
(mod 2) is a $k$-factor-critical graph. Zhai, Wei and Zhang \cite{ZWZ} presented some characterizations for the existence of $k$-factor-critical
graphs. More results on graph factors were obtained by Wang and Zhang \cite{WZo,WZs}, Zhou et al \cite{Zr1,ZZ,ZBS,ZPX1}, Gao, Wang and Chen \cite{GWC}.

\medskip

In this article, we also investigate the problem on the existence of $k$-factor-critical graphs, and characterize $k$-factor-critical graphs with
respect to the spectral radius. Our main result will be given in Sections 3.

\section{Preliminary Lemmas}

For a nonnegative integer $k$, a graph $G$ is said to be $k$-factor-critical if $G-Q$ admits a 1-factor for any $Q\subseteq V(G)$ with $|Q|=k$.

\medskip

\noindent{\textbf{Lemma 2.1}} (\cite{F}). Let $k$ be a nonnegative integer, and let $G$ be a graph of order $n$ with $n\equiv k$ (mod 2). Then
$G$ is $k$-factor-critical if and only if
$$
o(G-D)\leq|D|-k
$$
for any $D\subseteq V(G)$ with $|D|\geq k$.

\medskip

The following two results are very useful for the proof of the main theorem.

\medskip

\noindent{\textbf{Lemma 2.2}} (\cite{B}). Let $G$ be a connected graph, and let $H$ be a proper subgraph of $G$. Then $\rho(G)>\rho(H)$.

\medskip

Let $M$ be a real symmetric matrix whose rows and columns are indexed by $V=\{1,2,\cdots,n\}$. Suppose that $M$ can be written as
\begin{align*}
M=\left(
  \begin{array}{ccc}
    M_{11} & \cdots & M_{1s}\\
    \vdots & \ddots & \vdots\\
    M_{s1} & \cdots & M_{ss}\\
  \end{array}
\right)
\end{align*}
by means of partition $\pi: V=V_1\cup V_2\cup\cdots\cup V_s$ , wherein $M_{ij}$ is the submatrix (block) of $M$ derived by rows in $V_i$ and the
columns in $V_j$. We denote by $q_{ij}$ the average row sum of $M_{ij}$. Then matrix $M_{\pi}=(q_{ij})$ is said to be the quotient matrix of $M$.
If the row sum of every block $M_{ij}$ is a constant, then the partition is equitable.

\medskip

\noindent{\textbf{Lemma 2.3}} (\cite{YYSX}).  Let $M$ be a real matrix with an equitable partition $\pi$, and let $M_{\pi}$ be the corresponding
quotient matrix. Then every eigenvalue of $M_{\pi}$ is an eigenvalue of $M$. Furthermore, if $M$ is a nonnegative, then the largest eigenvalues of
$M$ and $M_{\pi}$ are equal.

\section{The main theorem and its proof}

In this section, we establish a relationship between spectral radius and $k$-factor-critical graphs.

\medskip

\noindent{\textbf{Theorem 3.1.}} Let $k$ be a nonnegative integer, and let $G$ be a $(k+1)$-connected graph of order $n$ with $n\equiv k$ (mod 2).\\
(\romannumeral1) For $n\geq k+4$ and $n\notin\{k+6,k+8\}$, or $(k,n)=(0,8)$, if $\rho(G)>\theta(n,k)$, then $G$ is $k$-factor-critical, where $\theta(n,k)$
is the largest root of $x^{3}-(n-4)x^{2}-(n+2k-1)x+2(k+1)(n-k-4)=0$.\\
(\romannumeral2) For $n=k+6$, if $\rho(G)>\frac{k+1+\sqrt{k^{2}+18k+33}}{2}$, then $G$ is $k$-factor-critical.\\
(\romannumeral3) For $k\geq1$ and $n=k+8$, if $\rho(G)>\frac{k+2+\sqrt{k^{2}+24k+64}}{2}$, then $G$ is $k$-factor-critical.

\medskip

\noindent{\it Proof.} Let $\varphi(x)=x^{3}-(n-4)x^{2}-(n+2k-1)x+2(k+1)(n-k-4)$ and let $\theta(n,k)$ be the largest root of $\varphi(x)=0$.
Suppose, to the contrary, that $G$ is not $k$-factor-critical. According to Lemma 2.1, there exists a subset $D\subseteq V(G)$ with $|D|\geq k$
such that $o(G-D)\geq|D|-k+1$. By parity, $o(G-D)\geq|D|-k+2\geq2$. Let $o(G-D)=\beta$ and $|D|=d$. Then we admit $\beta\geq d-k+2\geq2$.
Select a $(k+1)$-connected graph $G$ with $n$ vertices such that its spectral radius is as large as possible. In terms of Lemma 2.2 and the
choice of $G$, the induced subgraph $G[D]$ and all connected components of $G-D$ are complete graphs. Furthermore, $G=G[D]\vee(G-D)$.

\noindent{\bf Claim 1.} $d\geq k+1$.

\noindent{\it Proof.} Let $d=k$. Then $\omega(G-D)\geq o(G-D)\geq|D|-k+2=d-k+2=2$, which contradicts that $G$ is $(k+1)$-connected. This
completes the proof of Claim 1. \hfill $\Box$

\noindent{\bf Claim 2.} $G-D$ does not admit even components.

\noindent{\it Proof.} Assume that there exists an even component $G_e$ in $G-D$. Then we create a new graph $G^{(1)}$ by by joining $G_e$
and $G_o$ (that is, $G_e\vee G_o$), where $G_o$ is an odd component of $G-D$. Clearly, $G_e\vee G_o$ is an odd component of $G^{(1)}-D$ and
$o(G^{(1)}-D)=o(G-D)\geq|D|-k+2$. Furthermore, $G$ is a proper subgraph of $G^{(1)}$. In view of Lemma 2.2, $\rho(G)<\rho(G^{(1)})$, which
contradicts the choice of $G$. This completes the proof of Claim 2. \hfill $\Box$

Let $G_1,G_2,\cdots,G_{\beta}$ be the odd components in $G-D$ with $|V(G_1)|=n_1\geq|V(G_2)|=n_2\geq\cdots\geq|V(G_{\beta})|=n_{\beta}$.
By virtue of Claim 2, there exists the partition $\{D,V(G_1),V(G_2),\cdots,V(G_{\beta})\}$ of $G$. Thus, the quotient matrix of the
partition $\{D,V(G_1),V(G_2),\cdots,V(G_{\beta})\}$ of $G$ equals
\begin{align*}
\left(
  \begin{array}{ccccc}
    d-1 & n_1 & n_2 & \cdots & n_{\beta}\\
    d & n_1-1 & 0 & \cdots & 0\\
    d & 0 & n_2-1 & \cdots & 0\\
    \vdots & \vdots & \vdots & \vdots & \vdots\\
    d & 0 & 0 & \cdots & n_{\beta}-1\\
  \end{array}
\right).
\end{align*}
Then the characteristic polynomial of the matrix is equal to
\begin{align*}
f_1(x)&=(x-d+1)(x-n_1+1)\cdots(x-n_{\beta}+1)-dn_1(x-n_2+1)\cdots(x-n_{\beta}+1)\\
&+dn_2(x-n_1+1)(x-n_3+1)\cdots(x-n_{\beta}+1)+\cdots\\
&+(-1)^{i}dn_i(x-n_1+1)\cdots(x-n_{i-1}+1)(x-n_{i+1}+1)\cdots(x-n_{\beta}+1)+\cdots\\
&+(-1)^{\beta}dn_{\beta}(x-n_1+1)\cdots(x-n_{\beta-1}+1).
\end{align*}

Since the partition $\{D,V(G_1),V(G_2),\cdots,V(G_{\beta})\}$ is equitable, it follows from Lemma 2.3 that the largest root, say $\rho_1$,
of $f_1(x)=0$ equals the spectral radius of $G$. Thus, we possess $\rho(G)=\rho_1$. Since $K_{d+n_1}$ is a proper subgraph of $G$, it follows
from Lemma 2.2 that $\rho_1=\rho(G)>\rho(K_{d+n_1})=n_1+d-1$.

\noindent{\bf Claim 3.} $n_2=n_3=\cdots=n_{\beta}=1$.

\noindent{\it Proof.} We first verify $n_{\beta}=1$. Assume that $n_{\beta}\geq3$. If $n_1=1$, then we deduce $1=n_1\geq n_2\geq\cdots\geq n_{\beta}\geq3$,
which is a contradiction. Next, we deal with $n_1\geq3$.

We create a new graph $G^{(2)}$ by deleting two vertices in $G_{\beta}$ and adding two vertices to $G_1$ by joining the two vertices to the
vertices in $V(G_1)\cup D$. Obviously, $G[D]$ and all connected components in $G^{(2)}-D$ are complete graphs. In what follows, we prove
$\rho(G)<\rho(G^{(2)})$.

Assume that $V(G_1')$ and $V(G_{\beta}')$ are the vertex sets obtained from $V(G_1)$ by adding the two vertices and from $V(G_{\beta})$ by
deleting the two vertices, respectively. The quotient matrix of the partition $\{D,V(G_1'),V(G_2),\cdots,V(G_{\beta-1}),V(G_{\beta}')\}$ of
$G^{(2)}$ admits the characteristic polynomial $f_2(x)$ obtained from $f_1(x)$ by replacing $n_1$ and $n_{\beta}$ by $n_1+2$ and $n_{\beta}-2$,
respectively. Note that the partition $\{D,V(G_1'),V(G_2),\cdots,V(G_{\beta-1}),V(G_{\beta}')\}$ is equitable. It follows from Lemma 2.3 that
the largest root, say $\rho_2$, of $f_2(x)=0$ equals the spectral radius of $G^{(2)}$. Thus, we derive $\rho(G^{(2)})=\rho_2$. According to
$f_1(\rho_1)=0$, $n_1\geq n_{\beta}$ and $\rho_1>n_1+d-1$, we infer
\begin{align*}
f_2(\rho_1)&=2(\rho_1-d+1)(\rho_1-n_2+1)\cdots(\rho_1-n_{\beta-1}+1)(n_\beta-n_1-2)\\
&-2d(\rho_1-n_2+1)\cdots(\rho_1-n_{\beta-1}+1)(\rho_1+n_1-n_{\beta}+2)\\
&+2dn_2(\rho_1-n_3+1)\cdots(\rho_1-n_{\beta-1}+1)(n_{\beta}-n_1-2)\\
&-2dn_3(\rho_1-n_2+1)(\rho_1-n_4+1)\cdots(\rho_1-n_{\beta-1}+1)(n_{\beta}-n_1-2)+\cdots\\
&+(-1)^{\beta}2d(\rho_1-n_2+1)\cdots(\rho_1-n_{\beta-1}+1)(n_{\beta}-n_1-2)<0,
\end{align*}
which yields $\rho(G)=\rho_1<\rho_2=\rho(G^{(2)})$, which contradicts the choice of $G$. Hence, we infer $n_{\beta}=1$. Similarly, we may verify
$n_2=n_3=\cdots=n_{\beta-1}=1$. This completes the proof of Claim 3. \hfill $\Box$

\noindent{\bf Claim 4.} $\beta=d-k+2$.

\noindent{\it Proof.} Assume that $\beta\geq d-k+4$. Then we establish a new graph $G^{(3)}$ obtained from $G$ by adding an edge to join $G_{\beta-1}$
and $G_{\beta}$ ($G_{\beta-1}$ and $G_{\beta}$ are two vertices by Claim 3). It is obvious that $o(G^{(3)}-D)\geq d-k+2$ and $G$ is a proper subgraph
of $G^{(3)}$. Then it follows from Lemma 2.2 that $\rho(G)<\rho(G^{(3)})$, which contradicts the choice of $G$. Hence, we deduce $\beta\leq d-k+2$
by parity. On the other hand, $\beta\geq d-k+2$. Thus, we infer $\beta=d-k+2$. This completes the proof of Claim 4. \hfill $\Box$

In what follows, we consider two cases by the value of $n_1$.

\noindent{\bf Case 1.} $n_1=1$.

In this case, we possess $G=K_d\vee(n-d)K_1=K_d\vee(d-k+2)K_1$ and $n=d+\beta=d+d-k+2=2d-k+2$ by Claims 2--4. Note that the quotient matrix of the
adjacency matrix of $G$ with the partition $\{V(K_d),V((n-d)K_1)\}$ is equal to
\begin{align*}
\left(
  \begin{array}{cc}
    d-1 & n-d\\
    d & 0\\
  \end{array}
\right).
\end{align*}
Then the characteristic polynomial of the matrix is equal to
\begin{align*}
f_3(x)&=x(x-d+1)-d(n-d)\\
&=x^{2}-(d-1)x-d(n-d).
\end{align*}
Note that the partition $\{V(K_d),V((n-d)K_1)\}$ is equitable. By virtue of Lemma 2.3, the largest root, say $\rho_3$, of $f_3(x)=0$ is equal to the spectral
radius of $G$. Thus, we possess
\begin{align}\label{eq:3.1}
\rho(G)=\rho(K_d\vee(n-d)K_1)=\rho_3=\frac{d-1+\sqrt{(d-1)^{2}+4d(n-d)}}{2}.
\end{align}

Recall that $n=2d-k+2$. If $d=k+1$, then $n=k+4$ and
$\rho(G)=\rho_3=\frac{d-1+\sqrt{(d-1)^{2}+4d(n-d)}}{2}=\frac{k+\sqrt{k^{2}+12k+12}}{2}=\theta(k+4,k)$,
which contradicts $\rho(G)>\theta(n,k)$ for $n=k+4$. If $d=k+2$, then $n=k+6$ and $\rho(G)=\rho_3=\frac{d-1+\sqrt{(d-1)^{2}+4d(n-d)}}{2}=\frac{k+1+\sqrt{k^{2}+18k+33}}{2}$, which contradicts $\rho(G)>\frac{k+1+\sqrt{k^{2}+18k+33}}{2}$.
Next, we deal with $d\geq k+3$.

In light of \eqref{eq:3.1}, $f_3(\rho_3)=0$ and $n=2d-k+2$, we deduce
\begin{align}\label{eq:3.2}
\varphi(\rho_3)=&\varphi(\rho_3)-\rho_3f_3(\rho_3)\nonumber\\
=&\rho_3^{3}-(n-4)\rho_3^{2}-(n+2k-1)\rho_3+2(k+1)(n-k-4)\nonumber\\
&-\rho_3^{3}+(d-1)\rho_3^{2}+d(n-d)\rho_3\nonumber\\
=&-(n-d-3)\rho_3^{2}+((d-1)n-d^{2}-2k+1)\rho_3+2(k+1)(n-k-4)\nonumber\\
=&-(d-k-1)\rho_3^{2}+(d+1)(d-k-1)\rho_3+4(k+1)(d-k-1)\nonumber\\
=&(d-k-1)(-\rho_3^{2}+(d+1)\rho_3+4k+4).
\end{align}
Let $g_3(\rho_3)=-\rho_3^{2}+(d+1)\rho_3+4k+4$. Together with \eqref{eq:3.1} and $n=2d-k+2$, we have
\begin{align}\label{eq:3.3}
g_3(\rho_3)=&-\rho_3^{2}+(d+1)\rho_3+4k+4\nonumber\\
=&-\left(\frac{d-1+\sqrt{(d-1)^{2}+4d(n-d)}}{2}\right)^{2}\nonumber\\
&+(d+1)\cdot\frac{d-1+\sqrt{(d-1)^{2}+4d(n-d)}}{2}+4k+4\nonumber\\
=&-d^{2}-d+kd+4k+3+\sqrt{(d-1)^{2}+4d(d-k+2)}.
\end{align}

If $d=k+3$, then it follows from \eqref{eq:3.3} that
\begin{align*}
g_3(\rho_3)=&-d^{2}-d+kd+4k+3+\sqrt{(d-1)^{2}+4d(d-k+2)}\\
=&-9+\sqrt{k^{2}+24k+64}.
\end{align*}
For $d=k+3$ and $k=0$, we deduce $n=8$ and $g_3(\rho_3)<0$. Combining this with \eqref{eq:3.2}, we admit $\varphi(\rho_3)=(d-k-1)g_3(\rho_3)<0$, which
yields $\rho(G)=\rho_3<\theta(n,k)$ for $(n,k)=(8,0)$. For $d=k+3$ and $k\geq1$, we have $n=k+8$ and
$\rho(G)=\rho_3=\frac{d-1+\sqrt{(d-1)^{2}+4d(n-d)}}{2}=\frac{k+2+\sqrt{k^{2}+24k+64}}{2}$, which contradicts $\rho(G)>\frac{k+2+\sqrt{k^{2}+24k+64}}{2}$.

\noindent{\bf Claim 5.} If $d\geq k+4$, then $d^{2}+d-kd-4k-3>\sqrt{(d-1)^{2}+4d(d-k+2)}$.

\noindent{\it Proof.} By a direct computation, we obtain
\begin{align}\label{eq:3.4}
(d^{2}+d-&kd-4k-3)^{2}-((d-1)^{2}+4d(d-k+2))\nonumber\\
=&d^{4}+(2-2k)d^{3}+(k^{2}-10k-10)d^{2}+(8k^{2}+2k-12)d+16k^{2}+24k+8\nonumber\\
:=&h_3(d),
\end{align}
where $h_3(d)=d^{4}+(2-2k)d^{3}+(k^{2}-10k-10)d^{2}+(8k^{2}+2k-12)d+16k^{2}+24k+8$. Let
$h_3(x)=x^{4}+(2-2k)x^{3}+(k^{2}-10k-10)x^{2}+(8k^{2}+2k-12)x+16k^{2}+24k+8$ be a real function in $x$ with $x\in[k+4,+\infty)$. The derivative
function of $h_3(x)$ is
$$
h_3'(x)=4x^{3}+3(2-2k)x^{2}+2(k^{2}-10k-10)x+8k^{2}+2k-12.
$$
Furthermore, we possess
$$
h_3''(x)=12x^{2}+6(2-2k)x+2(k^{2}-10k-10).
$$
Note that
$$
-\frac{6(2-2k)}{24}=\frac{k-1}{2}<k+4\leq d.
$$
Hence, $h_3''(x)$ is increasing in the interval $[k+4,+\infty)$. Thus $h_3''(x)\geq h_3''(k+4)=2k^{2}+40k+220>0$, which implies that $h_3'(x)$ is
increasing in the interval $[k+4,+\infty)$ and so $h_3'(x)\geq h_3'(k+4)=2k^{2}+46k+260>0$. Thus, we infer that $h_3(x)$ is increasing in the interval
$[k+4,+\infty)$. Together with $d\geq k+4$, we derive
$$
h_3(d)\geq h_3(k+4)=4k+184>0.
$$
Combining this with \eqref{eq:3.4}, we deduce $(d^{2}+d-kd-4k-3)^{2}>(d-1)^{2}+4d(d-k+2)$, that is, $d^{2}+d-kd-4k-3>\sqrt{(d-1)^{2}+4d(d-k+2)}$.
This completes the proof of Claim 5. \hfill $\Box$

If $d\geq k+4$, then it follows from \eqref{eq:3.2}, \eqref{eq:3.3} and Claim 5 that
\begin{align*}
\varphi(\rho_3)=&(d-k-1)g_3(\rho_3)\\
=&(d-k-1)(-d^{2}-d+kd+4k+3+\sqrt{(d-1)^{2}+4d(d-k+2)})\\
<&0,
\end{align*}
which leads to $\rho(G)=\rho_3<\theta(n,k)$, which contradicts $\rho(G)>\theta(n,k)$.

\noindent{\bf Case 2.} $n_1\geq3$.

In this case, we infer $G=K_d\vee(K_{n_1}\cup(\beta-1)K_1)=K_d\vee(K_{n_1}\cup(d-k+1)K_1)$ and $n=d+n_1+d-k+1=2d+n_1-k+1$ by Claims 2--4. In terms
of the partition $\{V(K_d),V(K_{n_1}),V((d-k+1)K_1)\}$, the quotient matrix of the adjacency matrix of $G$ is equal to
\begin{align*}
\left(
  \begin{array}{ccc}
    d-1 & n-2d+k-1 & d-k+1\\
    d & n-2d+k-2 & 0\\
    d & 0 & 0\\
  \end{array}
\right).
\end{align*}
Then the characteristic polynomial of the matrix equals
\begin{align*}
f_4(x)=x^{3}-(n-d+k-3)x^{2}-(n+d^{2}-kd+k-2)x+d(d-k+1)(n-2d+k-2).
\end{align*}
Note that the partition $\{V(K_d),V(K_{n_1}),V((d-k+1)K_1)\}$ is equitable. According to Lemma 2.3, the largest root, say $\rho_4$, of $f_4(x)=0$
equals the spectral radius of $G$. Thus, we have $\rho(G)=\rho_4$ and $f_4(\rho_4)=0$.

Next we are to verify $\varphi(\rho_4)<0$. By plugging the value $\rho_4$ into $x$ of $\varphi(x)-f_4(x)$, we possess
\begin{align}\label{eq:3.5}
\varphi(\rho_4)=&\varphi(\rho_4)-f_4(\rho_4)\nonumber\\
=&(d-k-1)(-\rho_4^{2}+(d+1)\rho_4-(d+2)n+2d^{2}-(k-6)d+2k+8)\nonumber\\
=&(d-k-1)g_4(\rho_4),
\end{align}
where $g_4(\rho_4)=-\rho_4^{2}+(d+1)\rho_4-(d+2)n+2d^{2}-(k-6)d+2k+8$. Note that $d\geq k+1$ and $n=2d+n_1-k+1\geq2d-k+4$. Then we obtain
\begin{align}\label{eq:3.6}
\frac{d+1}{2}<\frac{d-1+\sqrt{(d-1)^{2}+4d(n-d)}}{2}.
\end{align}

Since $K_d\vee(n-d)K_1$ is a proper subgraph of $G$, we infer
\begin{align}\label{eq:3.7}
\rho_4=\rho(G)>\rho(K_d\vee(n-d)K_1)=\frac{d-1+\sqrt{(d-1)^{2}+4d(n-d)}}{2}
\end{align}
by \eqref{eq:3.1} and Lemma 2.2. Recall that $g_4(\rho_4)=-\rho_4^{2}+(d+1)\rho_4-(d+2)n+2d^{2}-(k-6)d+2k+8$. In terms of \eqref{eq:3.6} and
\eqref{eq:3.7}, we get
$$
\frac{d+1}{2}<\frac{d-1+\sqrt{(d-1)^{2}+4d(n-d)}}{2}<\rho(G)=\rho_4.
$$
Hence, we infer
\begin{align}\label{eq:3.8}
g_4(\rho_4)<&g_4\left(\frac{d-1+\sqrt{(d-1)^{2}+4d(n-d)}}{2}\right)\nonumber\\
=&-(2d+2)n+3d^{2}+(7-k)d+2k+7+\sqrt{(d-1)^{2}+4d(n-d)}.
\end{align}

\noindent{\bf Claim 6.} If $d\geq k+1$ and $n\geq2d-k+4$, then $(2d+2)n-3d^{2}-(7-k)d-2k-7>\sqrt{(d-1)^{2}+4d(n-d)}$.

\noindent{\it Proof.} By a direct computation, we possess
\begin{align}\label{eq:3.9}
((2d+2)n&-3d^{2}-(7-k)d-2k-7)^{2}-((d-1)^{2}+4d(n-d))\nonumber\\
=&(2d+2)^{2}n^{2}-(12d^{3}+(40-4k)d^{2}+(60+4k)d+8k+28)n\nonumber\\
&+9d^{4}+3(14-2k)d^{3}+(k^{2}-2k+94)d^{2}\nonumber\\
&+(-4k^{2}+14k+100)d+4k^{2}+28k+48\nonumber\\
:=&h_4(n),
\end{align}
where $h_4(n)=(2d+2)^{2}n^{2}-(12d^{3}+(40-4k)d^{2}+(60+4k)d+8k+28)n+9d^{4}+3(14-2k)d^{3}+(k^{2}-2k+94)d^{2}+(-4k^{2}+14k+100)d+4k^{2}+28k+48$.
Recall that $n\geq2d-k+4$. Then
$$
\frac{12d^{3}+(40-4k)d^{2}+(60+4k)d+8k+28}{2(2d+2)^{2}}<2d-k+4\leq n.
$$
Hence, we deduce
\begin{align}\label{eq:3.10}
h_4(n)\geq&h_4(2d-k+4)\nonumber\\
=&d^{4}+(10-2k)d^{3}+(k^{2}-18k+22)d^{2}+(8k^{2}-38k-4)d+16k^{2}-8k\nonumber\\
:=&l_4(d)
\end{align}
Let
$l_4(x)=x^{4}+(10-2k)x^{3}+(k^{2}-18k+22)x^{2}+(8k^{2}-38k-4)x+16k^{2}-8k$ be a real function in $x$ with $x\in[k+1,+\infty)$. The derivative
function of $l_4(x)$ is
$$
l_4'(x)=4x^{3}+3(10-2k)x^{2}+2(k^{2}-18k+22)x+8k^{2}-38k-4.
$$
Furthermore, we have
$$
l_4''(x)=12x^{2}+6(10-2k)x+2(k^{2}-18k+22).
$$
Note that
$$
-\frac{6(10-2k)}{24}=\frac{k-5}{2}<k+1\leq d.
$$
Hence, $l_4''(x)$ is increasing in the interval $[k+1,+\infty)$. Thus $l_4''(x)\geq l_4''(k+1)=2k^{2}+36k+116>0$, which implies that $l_4'(x)$ is
increasing in the interval $[k+1,+\infty)$ and so $l_4'(x)\geq l_4'(k+1)=4k^{2}+36k+74>0$. Thus, we infer that $l_4(x)$ is increasing in the interval
$[k+1,+\infty)$. Together with $d\geq k+1$, we obtain
$$
l_4(d)\geq l_4(k+1)=3k^{2}+8k+29>0.
$$
Combining this with \eqref{eq:3.9} and \eqref{eq:3.10}, we deduce $((2d+2)n-3d^{2}-(7-k)d-2k-7)^{2}>(d-1)^{2}+4d(n-d)$, that is,
$(2d+2)n-3d^{2}-(7-k)d-2k-7>\sqrt{(d-1)^{2}+4d(n-d)}$. This completes the proof of Claim 6. \hfill $\Box$

If $d\geq k+1$ and $n\geq2d-k+4$, then it follows from \eqref{eq:3.5}, \eqref{eq:3.8} and Claim 6 that
\begin{align*}
\varphi(\rho_4)=&(d-k-1)g_4(\rho_4)\\
\leq&(d-k-1)(-(2d+2)n+3d^{2}+(7-k)d+2k+7+\sqrt{(d-1)^{2}+4d(n-d)})\\
\leq&0,
\end{align*}
which yields $\rho(G)=\rho_4\leq\theta(n,k)$, which contradicts $\rho(G)>\theta(n,k)$ for $n\geq k+4$ and $n\notin\{k+6,k+8\}$, or $(k,n)=(0,8)$.

As for $n=k+6$, one has $\varphi(x)=x^{3}-(k+2)x^{2}-(3k+5)x+4(k+1)$ and $\varphi'(x)=3x^{2}-2(k+2)x-3k-5$. By a direct computation, we have
$\varphi\Big(\frac{k+1+\sqrt{k^{2}+18k+33}}{2}\Big)=k-3+\sqrt{k^{2}+18k+33}>0$ and
$\varphi'\Big(\frac{k+1+\sqrt{k^{2}+18k+33}}{2}\Big)=\frac{k^{2}+18k+37+(k-1)\sqrt{k^{2}+18k+33}}{2}>0$, and so
$\rho(G)=\rho_4\leq\theta(n,k)<\frac{k+1+\sqrt{k^{2}+18k+33}}{2}$, which is a contradiction to $\rho(G)>\frac{k+1+\sqrt{k^{2}+18k+33}}{2}$ for
$n=k+6$.

As for $k\geq1$ and $n=k+8$, one has $\varphi(x)=x^{3}-(k+4)x^{2}-(3k+7)x+8(k+1)$ and $\varphi'(x)=3x^{2}-2(k+4)x-3k-7$. By a direct calculation,
we derive $\varphi\Big(\frac{k+2+\sqrt{k^{2}+24k+64}}{2}\Big)=-18+2\sqrt{k^{2}+24k+64}>0$ and
$\varphi'\Big(\frac{k+2+\sqrt{k^{2}+24k+64}}{2}\Big)=\frac{k^{2}+24k+72+(k-2)\sqrt{k^{2}+24k+64}}{2}>0$, and so
$\rho(G)=\rho_4\leq\theta(n,k)<\frac{k+2+\sqrt{k^{2}+24k+64}}{2}$, which is a contradiction to $\rho(G)>\frac{k+2+\sqrt{k^{2}+24k+64}}{2}$ for
$k\geq1$ and $n=k+8$. This completes the proof of Theorem 3.1. \hfill $\Box$

\section{Extremal graphs}

In this section, we claim that the spectral radius conditions in Theorem 3.1 are sharp.

\medskip

\noindent{\textbf{Theorem 4.1.}} Let $k$ and $n$ be two nonnegative integers with $n\equiv k$ (mod 2), and let $\theta(n,k)$ be the largest root
of $x^{3}-(n-4)x^{2}-(n+2k-1)x+2(k+1)(n-k-4)=0$. For $n\geq k+4$ and $n\notin\{k+6,k+8\}$, or $(k,n)=(0,8)$, we possess
$\rho(K_{n-k-3}\vee K_{k+1}\vee(2K_1))=\theta(n,k)$ and $K_{n-k-3}\vee K_{k+1}\vee(2K_1)$ is not a $k$-factor-critical graph. For $n=k+6$, we
admit $\rho(K_{k+2}\vee(4K_1))=\frac{k+1+\sqrt{k^{2}+18k+33}}{2}$ and $K_{k+2}\vee(4K_1)$ is not a $k$-factor-critical graph. For $n=k+8$, we
have $\rho(K_{k+3}\vee(5K_1))=\frac{k+2+\sqrt{k^{2}+24k+64}}{2}$ and $K_{k+3}\vee(5K_1)$ is not a $k$-factor-critical graph.
\medskip

\noindent{\it Proof.} Consider the vertex partition $\{V(K_{n-k-3}),V(K_{k+1}),V(2K_1)\}$ of $K_{n-k-3}\vee K_{k+1}\vee(2K_1)$. The quotient matrix
of the partition $\{V(K_{n-k-3}),V(K_{k+1}),V(2K_1)\}$ of $K_{n-k-3}\vee K_{k+1}\vee(2K_1)$ is equal to
\begin{align*}
\left(
  \begin{array}{ccc}
    n-k-4 & k+1 & 0\\
    n-k-3 & k & 2\\
    0 & k+1 & 0\\
  \end{array}
\right).
\end{align*}
Then the characteristic polynomial of the matrix equals $x^{3}-(n-4)x^{2}-(n+2k-1)x+2(k+1)(n-k-4)$. Note that the partition is equitable. Then it
follows from Lemma 2.3 that the largest root of $x^{3}-(n-4)x^{2}-(n+2k-1)x+2(k+1)(n-k-4)=0$ equals the spectral radius of the graph
$K_{n-k-3}\vee K_{k+1}\vee(2K_1)$. That is, $\rho(K_{n-k-3}\vee K_{k+1}\vee(2K_1))=\theta(n,k)$. Set $D=V(K_{k+1})$. Then
$o(K_{n-k-3}\vee K_{k+1}\vee(2K_1)-D)=3>1=(k+1)-k=|D|-k$. According to Lemma 2.1, $K_{n-k-3}\vee K_{k+1}\vee(2K_1)$ is not a $k$-factor-critical
graph.

\medskip

Consider the vertex partition $\{V(K_{k+2}),V(4K_1)\}$ of $K_{k+2}\vee(4K_1)$. The quotient matrix of the partition $\{V(K_{k+2}),V(4K_1)\}$ of
$K_{k+2}\vee(4K_1)$ equals
\begin{align*}
\left(
  \begin{array}{cc}
    k+1 & 4\\
    k+2 & 0\\
  \end{array}
\right).
\end{align*}
Then the characteristic polynomial of the matrix is $x^{2}-(k+1)x-4(k+2)$. Note that the partition is equitable. In view of Lemma 2.3, the largest
root of $x^{2}-(k+1)x-4(k+2)=0$ is equal to the spectral radius of the graph $K_{k+2}\vee(4K_1)$. Thus, we derive
$\rho(K_{k+2}\vee(4K_1))=\frac{k+1+\sqrt{k^{2}+18k+33}}{2}$. Set $D=V(K_{k+2})$. Then $o(K_{k+2}\vee(4K_1)-D)=4>2=(k+2)-k=|D|-k$. By means of Lemma
2.1, $K_{k+2}\vee(4K_1)$ is not a $k$-factor-critical graph.

\medskip

Consider the vertex partition $\{V(K_{k+3}),V(5K_1)\}$ of $K_{k+3}\vee(5K_1)$. The quotient matrix of the partition $\{V(K_{k+3}),V(5K_1)\}$ of
$K_{k+3}\vee(5K_1)$ is equal to
\begin{align*}
\left(
  \begin{array}{cc}
    k+2 & 5\\
    k+3 & 0\\
  \end{array}
\right).
\end{align*}
Then the characteristic polynomial of the matrix equals $x^{2}-(k+2)x-5(k+3)$. Note that the partition is equitable. According to Lemma 2.3, the
largest root of $x^{2}-(k+2)x-5(k+3)=0$ equals the spectral radius of the graph $K_{k+3}\vee(5K_1)$. Thus, we get
$\rho(K_{k+3}\vee(5K_1))=\frac{k+2+\sqrt{k^{2}+24k+64}}{2}$. Set $D=V(K_{k+3})$. Then $o(K_{k+3}\vee(5K_1)-D)=5>3=(k+3)-k=|D|-k$. In terms of Lemma
2.1, $K_{k+3}\vee(5K_1)$ is not a $k$-factor-critical graph. \hfill $\Box$

\medskip

\section*{Data availability statement}

My manuscript has no associated data.

\section*{Declaration of competing interest}

The authors declare that they have no conflicts of interest to this work.


\end{document}